\documentclass[20pt, twoside, a4paper]{article}

\let\Finishall\relax            
\ifx\TestIngCommand\undefined\relax                                     
\else\input ../../../proc/ppreamb-book.tex                               
\input init.tex \fi                                                     
\let\TestIngCommand\undefined                                             

\newtheorem{example}{Example}
\newtheorem{question}{Question}
\newtheorem{remark}{Remark}

\usepackage{amsmath,amscd,amssymb}
\usepackage{graphics}
\newtheorem{theo}{Theorem}

\newtheorem{prop}{Proposition}
\newtheorem{defi}{Definition}

\newskip\ttglue\ttglue=.5em plus.25em minus.15em
\def\firstname#1{\def\FIRSTNAME{#1}\ignorespaces}                  
\def\lastname#1{\def\LASTNAME{#1}\ignorespaces}
\def\middleinitial#1{\def\MIDDLEINI{#1}\ignorespaces}
\def\department#1{\def\DEPARTMENT{#1}\ignorespaces}
\def\institute#1{\def\INSTITUTE{#1}\ignorespaces}
\def\address#1{\def\ADDRESS{#1}\ignorespaces}
\def\country#1{\def\COUNTRY{#1}\ignorespaces}
\def\otheraffiliation#1{\def\OTHERAFFILIATION{#1}\ignorespaces}
\def\email#1{\def\EMAIL{#1}\ignorespaces}
\newcount\autcount\autcount=0
\newcount\affcount\affcount=0
\newcount\numcount\newcount\nummcount
\newcount\nummmcount\newcount\nummmmcount
\def\writename#1#2{\ \kern-1ex\hbox{
  \csname AUthor\the#1\endcsname\
  \edef\TESTSTR{}\expandafter\ifx\csname auTHor\the#1\endcsname\TESTSTR
  \else\csname auTHor\the#1\endcsname.\ \fi 
  \csname authOR\the#1\endcsname$^{\csname AFF\the#1\endcsname}$
  \expandafter\ifx\csname corr\number#1\endcsname\relax
  \else\thanks{Corresponding author.}\ \fi
  }\ifnum#1<#2, \else\ \kern-1ex\fi}
\def\writeemail#1{
  \nummcount=0\relax\nummmcount=0\relax
  \loop\ifnum\nummcount<\autcount\advance\nummcount by1\relax
    {\expandafter\ifnum\csname AFF\the\nummcount\endcsname=#1\relax
    \global\advance\nummmcount by1\fi}\repeat
  \nummcount=0\relax\nummmmcount=0\relax
  \loop\ifnum\nummcount<\autcount\advance\nummcount by1\relax
    {\expandafter\ifnum\csname AFF\the\nummcount\endcsname=#1\relax
    \global\advance\nummmmcount by1\relax\def\blank{}\expandafter
    \ifx\csname EMAIL\the\nummcount\endcsname\blank(no e-mail)
    \else\csname EMAIL\the\nummcount\endcsname
    \fi 
    \ifnum\nummmmcount<\nummmcount; \fi\fi}\repeat}
\long\def\BeginAuthorList#1\EndAuthorList{#1\relax
  \author{\vbox{\hsize=390pt\noindent\numcount=0\relax
    \loop\ifnum\numcount<\autcount\advance\numcount by1\relax
      \writename{\numcount}{\autcount}
      \repeat}\\[2mm]
    \vbox{\small\numcount=0\relax
      \loop\ifnum\numcount<\affcount\advance\numcount by1\relax
        \vbox{{\count0=\numcount\relax
          \loop\expandafter\ifnum\csname AFF\the\count0\endcsname
            <\numcount\relax\advance\count0 by1\relax\repeat
          $^{\csname AFF\the\count0\endcsname}$}
        \def\BLANK{}\expandafter\ifx\csname DEPT\the\numcount\endcsname
          \BLANK
          \else\csname DEPT\the\numcount\endcsname, \fi
        \csname INST\the\numcount\endcsname,
        \csname ADDR\the\numcount\endcsname,
        \csname COUN\the\numcount\endcsname
        \edef\TEST{}\expandafter\ifx\csname OTHE\the\numcount\endcsname
          \TEST
          .\else;\break\csname OTHE\the\numcount\endcsname.\fi}
        \vbox{\writeemail{\numcount}}
        \repeat}\\}}
\expandafter\def\csname x1\endcsname{}
\expandafter\def\csname x2\endcsname{}
\expandafter\def\csname x3\endcsname{}
\expandafter\def\csname x4\endcsname{}
\expandafter\def\csname x5\endcsname{}
\expandafter\def\csname x6\endcsname{}
\expandafter\def\csname x7\endcsname{}
\expandafter\def\csname x8\endcsname{}
\expandafter\def\csname x9\endcsname{}
\def\Author#1#2{\global\advance\autcount by1\relax#2
  \expandafter\edef\csname AUthor\the\autcount\endcsname{\FIRSTNAME}
  \expandafter\edef\csname auTHor\the\autcount\endcsname{\MIDDLEINI}
  \expandafter\edef\csname authOR\the\autcount\endcsname{\LASTNAME}
  \expandafter\edef\csname EMAIL\the\autcount\endcsname{\EMAIL}
  \let\tempera\"\def\"{\string\"}\expandafter\ifx\csname x\DEPARTMENT
    \endcsname\relax
    \global\advance\affcount by1\relax\let\"\tempera
    \expandafter\edef\csname DEPT\the\affcount\endcsname{\DEPARTMENT}
    \expandafter\edef\csname INST\the\affcount\endcsname{\INSTITUTE}
    \expandafter\edef\csname ADDR\the\affcount\endcsname{\ADDRESS}
    \expandafter\edef\csname COUN\the\affcount\endcsname{\COUNTRY}
    \expandafter\edef\csname OTHE\the\affcount\endcsname{\OTHERAFFILIATION}
    \expandafter\edef\csname AFF\the\autcount\endcsname{\the\affcount}
  \else\expandafter\edef\csname AFF\the\autcount\endcsname{\DEPARTMENT}
  \fi\let\"\tempera\ignorespaces}
\def\CorrespondingAuthor#1#2{
  \expandafter\xdef\csname corr\number#1\endcsname{cor}
  \Author#1{#2}}
\def\PaperTitle#1{\title{\bf#1}}
\def\Category#1{\ignorespaces}
\def\keywords#1{{\noindent \emph{Keywords:}
  \def\BLANK{}\def\TEST{#1}\ifx\BLANK\TEST(n/a).\else#1\fi}}
\begin{document}
\begin{LARGE}
\large{                                                          
\PaperTitle{Multidimensional Fourier Quasicrystals I. Sufficient Conditions}

\Category{(Pure) Mathematics}

\date{}

\BeginAuthorList
\Author1{
\firstname{Wayne}
    \lastname{Lawton}
    \middleinitial{M}
    \department{Department of the Theory of Functions, Institute of Mathematics and Computer Science}
    \institute{Siberian Federal University}
    \otheraffiliation{}
   \address{Krasnoyarsk}
    \country{Russian Federation}
    \email{wlawton@gmail.com}}
    \Author2{
     \firstname{August}
    \lastname{Tsikh}
    \middleinitial{K}
    \department{Department of the Theory of Functions, Institute of Mathematics and Computer Science}
    \institute{Siberian Federal University}
    \otheraffiliation{}
   \address{Krasnoyarsk}
    \country{Russian Federation}
    \email{atsikh@sfu-kras.ru}}
\EndAuthorList
\maketitle
\thispagestyle{empty}
\begin{abstract}
We derive sufficient conditions for an atomic measure $\sum_{\lambda \in \Lambda} m_\lambda\, \delta_\lambda,$ where $\Lambda \subset \mathbb R^n,$ $m_\lambda$ are positive integers, and $\delta_\lambda$ is the point measure at $\lambda,$ to be a Fourier quasicrystal, and suggest why they may also be necessary. These conditions extend the necessary and sufficient conditions derived by Lev, Olevskii, and Ulanovskii for $n = 1.$ Our methods exploit the toric geometry relation between Grothendieck residues and Newton polytopes derived by Gelfond and Khovanskii. 
\end{abstract}
\noindent{\bf 2020 Mathematics Classification 52C23; 32A60; 32A27}
%
\footnote{\thanks{This work is supported by the Krasnoyarsk Mathematical 
Center and financed by the Ministry of Science and Higher Education 
of the Russian Federation in the framework of the establishment and 
development of Regional Centers for Mathematics Research and 
Education (Agreement No. 075-02-2020-1534/1).}}
%
%
\section{Introduction}\label{sec0}
A quasicrystal is an atomic measure 
\begin{equation}\label{eqn0.1}
\mu = \sum_{\lambda \in \Lambda} c_\lambda \delta_\lambda,
\end{equation} 
where $\Lambda \subset \mathbb R^n$ is discrete, $c_\lambda \in \mathbb C \backslash \{0\},$ and $\delta_\lambda$ is the point measure at $\lambda,$ that is
a tempered distribution whose Fourier transform is an atomic measure
\begin{equation}\label{eqn0.2} 
\widehat \mu = \sum_{s \in S} a_s \delta_s, 
\end{equation} 
where $a_s \in \mathbb C\backslash \{0\}.$ $\Lambda$ is the support of $\mu$ and $S \subset \mathbb R^n$ is the spectrum of $\mu.$ 
Alloys synthesized by Shechtman \cite{schechtman} and found in meteorites by Steinhardt \cite{steinhardt} exhibit discrete diffraction patterns with icosohedral symmetry impossible for conventional periodic crystals, so were named quasicrystals. Sets associated with physical quasicrystals and with Penrose and Ammann aperiodic planar tilings \cite{debruijn, grunbaumshephard} are model (cut and project) sets constructed earlier by Meyer 
\cite{meyer1, meyer2}. Each $\Lambda$ is uniformly discrete or Delone \cite{delone}, i.e. the distance between any two of its points is bounded below by a positive number, and $\mu := \sum_{\lambda \in \Lambda} \delta_\lambda$ is a quasicrystal with a dense spectrum. Furthermore, for every compactly supported continuous $f : \mathbb R^n \mapsto \mathbb C$ the convolution
\begin{equation}\label{eqn0.3} 
(\mu*f)(x) := \int_{\mathbb R^n} f(x-t)d\mu(t) 
\end{equation} 
is Besicovitch but not Bohr almost periodic \cite{besicovitch, bohr, favorov3}. Moody surveys model sets in \cite{moody}.
\\ \\
A Crystalline measure (CM) is a quasicrystal whose spectrum is discrete. A set is called relatively uniformly discrete if it is the union of a finite number of uniformly discrete sets. Poisson's summation formula implies that 
$\mu = \sum_{k \in \mathbb Z^n} \delta_k$ satisfies 
$\widehat \mu = \mu$ so $\mu$ is a CM. 
A simple computation then gives: if $J \geq 1$ and for $1 \leq j \leq J,$ 
$\tau_j \in \mathbb R^n,$ $L_j \subset \mathbb R^n$ is a 
rank--$n$ lattice subgroup, and
$f_j : \mathbb R^n \mapsto \mathbb C$ is a trigonometric 
polynomial, i.e. a finite linear combination of
exponentials $e^{2\pi i <\omega,x>},$ then
\begin{equation}\label{eqn0.4} 
	\mu := \sum_{j=1}^J \sum_{\lambda \in L_j+\tau_j} f_j(\lambda)\, \delta_\lambda
\end{equation} 
is a CM with both $\Lambda$ and $S$ relatively uniformly discrete. We call measures of the form (\ref{eqn0.4}) trivial. If all $L_j$ are equal then $\Lambda$ and $S$ are uniformly discrete. Lagarias (\cite{lagarias}, Problem 4.1) conjectured that if 
$\Lambda$ and $S$ are uniformly discrete then $\mu$ is trivial and all $L_j$ are equal. Lev and Olevskii \cite{levolevskii1} proved this for $n = 1$ or $\mu$ positive. Favorov \cite{favorov4} constructed a trivial measure for $n = 2,$ $J = 2$ and $L_1 \neq L_2$ for which both $\Lambda$ and $S$ are uniformly discrete and in \cite{favorov5} he gave sufficient conditions that imply all $L_j$ are equal. We conjecture that $\Lambda$ and $S$ are relatively uniformly discrete iff $\mu$ is trivial. Supports of trivial measures are unions of infinite arithmetic progressions, i.e. translates of a subroup isomorphic to $\mathbb Z.$ 
In \cite{levolevskii3} Lev and Olevskii
constructed a CM on $\mathbb R$ whose support does not contain any
infinite arithmetic 
progression, however both its support and spectrum generate finite dimensional subspaces over $\mathbb Q.$  Kolountzakis \cite{kolountzakis} constructing a CM whose support and spectrum generate infinite dimensional subspaces over $\mathbb Q.$ In the 1950's Weil \cite{weil}, Kahane and Mandelbrojt 
\cite{kahanemandelbrojt} and Guinand \cite{guinand} studied measures related to the Riemann hypothesis. 
In \cite{meyer3} Meyer used Guinand's measures to construct a nontrivial 
CM on $\mathbb R$ including one for which both
$\Lambda \cap (0,\infty)$ and $S \cap (0,\infty)$ are 
linearly independent over $\mathbb Q.$ Ronkin \cite{ronkin}
proved that every Bohr almost periodic measure is translation
bounded, i.e. its variation $|\mu|$ is uniformly bounded on balls of radius $1.$ Meyer's CMs are not translation bounded and hence not Bohr almost periodic.
\\ \\
A Fourier quasicrystal (FQ) is a CM $\mu$ whose variation of $\widehat \mu$
\begin{equation}\label{eqn0.5} 
|\widehat \mu| := \sum_{s \in S} |a_s| \delta_s
\end{equation}
is a tempered distribution, i.e. there exists $N > 0$ and $C > 0$ satisfying
\begin{equation}\label{eqn0.8} 
\sum_{s \in S \cap B(0,R)} |a_s| < C(1+R)^N, \ R > 0.
\end{equation}
The CM examples above are FQs but recently Favorov \cite{favorov6} constructed a CM that is not a FQ.  Olevski and Ulanovskii (\cite{olevskiiulanovskii2}, Proposition 4) proved that every positive FQ on $\mathbb R$ is translation bounded and their proof extends directly to $\mathbb R^n.$ 
Therefore (\cite{favorov5}, Theorem 11) implies that every positive FQ is a Bohr almost periodic measure. We define a multivariate Laurent polynomial $P(z_1,..,z_n)$ to be stable if it has no zeros whenever all $|z_j| < 1$ and pair--stable if in addition $P(z_1^{-1},...,z_n^{-1})$ is stable (different from Wagner's definition \cite{wagner} but related by a Cayley transformation).
Lee-Yang polynomials (\cite{leeyang}, Theorem 5.12), (\cite{passaretsikh}, Theorem 2) and quantum graph polynomials \cite{barragaspard,colin} are pair-stable. If $P$ is pair-stable and $\omega_1,...,\omega_n > 0$ then the univariate trigonometric polynomial $f(z) := P(e^{2\pi i\omega_1 z},...,e^{2\pi i\omega_n z})$ has only real roots.
Kurashov and Sarnak \cite{kurasovsarnak} proved that if $f$ is as above, $\Lambda$ is its zero set, and $m_\lambda$ is the multiplicity of $\lambda,$ then $\mu := \sum_{\lambda \in \Lambda} m_\lambda \, \delta_\lambda$ is a FQ. They constructed an example where the pair-stable polynomial
\begin{equation}\label{eqn0.9} 
P(z_1,z_2) := 1 - \frac{1}{3}z_1+\frac{1}{3}z_2^2 - z_1z_2^2, 
\end{equation}
$\omega_1, \omega_2 > 0,$ and $\omega_2/\omega_1$ irrational. They used a deep conjecture of Lang \cite{lang}, proved by Liardet \cite{liardet} for dimension $2$ and Laurent \cite{laurent} for dimension $\geq 3,$ to prove that $\Lambda$ generates an infinite dimensional vector space over $\mathbb Q.$ Evertse (\cite{evertse}, Theorem 10.10.1) gives a detailed explanation of Lang's conjecture and its proof. 
Olevskii and Ulanovskii (\cite{olevskiiulanovskii1}, Example 1) constructed a FQ  on $\mathbb R$ that is a special case of the following construction in (\cite{lawton}, Example 2): $\delta \in (-1,1)\backslash \{0\},$
\begin{equation}\label{eqn0.9} 
P(z_1,z_2) := z_1 - z_1^{-1} + \delta(z_2 - z_2^{-1}) .
\end{equation}
$f(z) := P(e^{2\pi i \omega_1\, z },e^{2\pi i \omega_2\, z}.$ 
In contrast to Kurasov and Sarnak's example, $P$ is not stable. Nevertheless $\Lambda \subset \mathbb R$ iff $|\omega_2/\omega_1| \leq 1$ and then 
$\mu := \sum_{\lambda \in \Lambda} m_\lambda \, \delta_\lambda$ is a FQ.
Furthermore, if $\omega_2/\omega_1$ is irrational then $\Lambda$ does not contain any infinite arithmetic progression. 
In (\cite{olevskiiulanovskii2}, Theorem 8)  Olevskii and Ulanovskii proved that for a 
measure of the form
$\sum_{\lambda \in \Lambda} m_\lambda \, \delta_\lambda,$ where $\Lambda \subset \mathbb R$ is discrete and $m_\lambda$ are positive integers, it is both necessary and sufficient that $\Lambda$ is the zero set of a trigonometric polynomial and $m_\lambda$ is the multiplicity of $\lambda.$ Sufficiency was proved in (\cite{olevskiiulanovskii1}, Corollary 1 of Theorem 1) based on earlier work of Lev and Olevski \cite{levolevskii2} which derived a generalized Poisson summation formula (PSF) for the zero set (not necessarily real) of a trigonometric polynomial. Their derivation computes a contour integral in the plane in two ways: directly and using Cauchy residues. We extend their sufficieny condition for $n \geq 2.$
\\ \\
Section \ref{sec1} introduces notation, derives representations of systems of trigonometric polynomials by Laurent polynomials, and relates zeros of Laurent polynomials to their Newton polytopes.  
Section \ref{sec2} gives Kazarnovskii's formula for the density of zeros of a system of $n$ trigomometric polynomials. 
Section \ref{sec3} derives a multidimensional generalized PSF. It uses two steps: (i) approximate by a periodic system of trigonometric polynomials represented by a system of Laurent polynomials, (ii) replace Cauchy residues by Grothendieck residues and apply the Gelfond--Khovanskii formula \cite{gelfondkhovanskii1,gelfondkhovanskii2} that equates sums of Grothendieck residues by weighted sums of residues at vertices of a Newton polytope.
Section \ref{sec4} uses the generalized PSF to derive a sufficient condition for a multidimensional FQ. It constructs two classes of multidimensional FQs, one similar to those construced by Meyer \cite{meyer4} using Ahern measures \cite{ahern} and the other similar to the one dimensional ones constructed by Olevskii and Ulanovslii.  It uses Lang's conjecture to prove that, unlike Meyer's two-dimensional FQ \cite{meyer5}, the supports of both classes of FQs do not contain any infinite arithmetic progression. 
Section \ref{sec5} formulates two quesions for future research. Question 1 concerns conditions under which the zeros of a system of trigonmetric polynomials are a subset of $\mathbb R^n.$ Question 2 suggest why our sufficient conditions for FQs may also be necessary.
%
%
\section{Preliminary Results}\label{sec1}}
$:=$ means is defined to equal, $\circ$ means composition of functions, and
$\simeq$ means isomorphism of topological groups.
$\mathbb N := \{1, 2, ...\}, \mathbb Z, \mathbb Q, \mathbb R, 
\mathbb C, \mathbb C^*$ are the natural, integer, rational, real, complex, and nonzero complex numbers. $\mathbb T := \{z \in \mathbb C:|z| = 1\}$  is the circle group. For $m, n \in \mathbb N,$
$\mathbb C^n$ is the complex vector space of column vectors,
$\mathbb R^n$ is its real subspace, and 
$\mathbb R^{m \times n}$ is the set of real $m$ by $n$ matrices
identified with $\mathbb R$--linear maps 
$M : \mathbb C^n \mapsto \mathbb C^m.$
If $a,b \in \mathbb R^n, [a,b] := \{ta+(1-t)b: t \in [0,1]\}$ is the line segment with endpoints $a$ and $b.$
$B(0,R) \subset \mathbb C^n$ is the open ball of radius $R$ centered at $0,$
$I_m \in \mathbb R^{m \times m}$ is the identity matrix,
$GL(n,\mathbb R)$ is the general linear group,
$GL(m,\mathbb Z) := \{U \in \mathbb Z^{m \times m}: \det U = \pm 1\}$ is the unimodular group. For $M \in \mathbb R^{m \times n}$ define
\begin{equation}\label{eqn1.1}
r(M) := \hbox{rank } \mathbb Z^m \cap M\mathbb R^n. 
\end{equation}
Clearly $r(M) \leq $rank$\, M.$ The subspace
$M\mathbb R^n$ is called rational if $r(M) = $rank$\,M.$ 
Define the operator norm
\begin{equation}\label{eqn1.2}
||A||_2 := \max\, \{\, ||Au||_2: u \in \mathbb R^m,\, ||u||_2 = 1\, \}, \ \ A \in \mathbb R^{m \times m}.
\end{equation}
$\mathbb C^{*m} := (\mathbb C^*)^m$ is the multiplicative torus group with identity $1,$ so $\mathbb T^m$ is its compact subgroup. 
$\rho_m : \mathbb C^m \rightarrow \mathbb C^{*m}$ is the epimorphism defined by
$
\rho_m(z)_j = e^{2\pi i z_j}, \, \, z \in \mathbb C^m, \, j =1,...,m.
$
For $N \in \mathbb Z^{m \times n},$ the homomorphism
$\widetilde N : \mathbb C^{*n} \mapsto \mathbb C^{*m}$ 
defined by
$(\widetilde N z)_i = \prod_{j=1}^n z_j^{N_{i,j}}, \, i = 1,...,m,$
satisfies
$\rho_m \circ N = \widetilde N \circ \rho_n \, : 
\, \mathbb C^n \mapsto \mathbb C^{*m}.$
\\ \\
$\mathcal T_n$ is the set of trigonometric polynomials 
$F = (f_1,...,f_n): \mathbb C^n \mapsto \mathbb C^n$
where each $f_j(z)$ is a linear combination, with nonzero complex coefficients, of $e^{2\pi i x\cdot z}$ where $x$ belongs to the spectrum $\Omega(f_j) \subset \mathbb R^n.$ 
\\
$\Gamma(F) := \{y \in \mathbb R^n: F(z+y) = F(z), \, z \in \mathbb C^n\}$ is the period group of $F,$
\\
$\Lambda(F) := \{z \in \mathbb C^n: F(z) = 0\}$ is the zero set of $F,$
\\
$G(F)$ is the subroup of  $\mathbb R^n$ generated by $\Omega(f_1) \cup \cdots \cup \Omega(f_n).$ 
\\
The Newton polytope $\mathcal N(f_j)$ of $f_j$ is the convex hull of $\Omega(f_j),$ $\mathcal V(f_j)$ its set of its vertices, and
$\mathcal N(F) := (\mathcal N(f_1),...,\mathcal N(f_n))$
\\ \\
$\mathcal L_{m,n}$ is the set of Laurent polynomials  
$P  = (p_1,...,p_n) :\mathbb C^{*m} \rightarrow \mathbb C^n$ 
where where each $p_j(z)$ is a linear combination, with nonzero complex coefficients, of 
$z^k := z_1^{k_1}\cdots z_m^{k_m}$ where $k$ belongs to the spectrum $\Omega(p_j) \subset \mathbb Z^m.$ $\mathcal N(p_j)$ is the Newton polytope of $p_j,$ $\mathcal V(p_j)$ is its set of vertices, and
$\mathcal N(P) = (\mathcal N(p_1),...,\mathcal N(p_n)).$ 
$\Lambda(P)$ is the zero set of $P.$ 
\begin{defi}\label{def1.1}
$\mathcal N(F)$ is unfolded if for every $y \in \mathbb R^n$ there exist $j \in \{1,...,,n\}$ and $v \in \mathcal V(f_j)$ satisfying the inequality
\begin{equation}\label{eqn1.3}
	y \cdot v < y \cdot u, \ \ u \in \mathcal V(f_j) \backslash \{v\}.
\end{equation}
\end{defi}
We record the following observations without proof.
\begin{prop}\label{prop1.1}
If $F \in \mathcal T_n,$ $Q \in \mathcal L_{n,n},$ and $B \in GL(n,\mathbb R),$ then
\begin{enumerate}\label{obs1.1}
\item if $\mathcal N(F)$ is unfolded it remains unfolded under small perturbations,
\item if $\mathcal N(F)$ is unfolded then so is
$\mathcal N(F \circ B) = (B^T\mathcal N(f_1),...,B^T\mathcal N(f_n)),$ 
\item if $\mathcal N(F)$ is unfolded, then $\Im \Lambda(F)$ is a bounded subset of $\mathbb R^n,$
\item $Q \circ \rho_n \in \mathcal T_n$ and $\mathcal N(Q \circ \rho_n) = \mathcal N(Q).$
\item if $\mathcal N(Q)$ is unfolded, then the analytic set 
$\Lambda(Q)$ is compact hence finite by (\cite{chirka}, Proposition 1, p. 31). 
\end{enumerate}
\end{prop}
The following result represents trigonometric by Laurent polynomials.
\begin{prop}\label{prop1.2}
\begin{enumerate}
\item Every continuous homomorphism $\psi : \mathbb R^n  \mapsto \mathbb T^m$ has the form $\psi = \rho_m \circ M$ where  $M \in \mathbb R^{m \times n}.$ The image $\psi(\mathbb R^n)$ is dense in $\mathbb T^m$ iff the rows of $M$ are linearly independent over $\mathbb Q.$ 
\item Ker$(\psi) \simeq \mathbb R^{\, n - \hbox{rank}\,M} \times \mathbb Z^{\, r(M)},$
hence  $\psi(\mathbb R^n) \simeq 
\mathbb R^{\, \hbox{rank}\,M - r(M)} \times \mathbb T^{\,r(M)}$
and $\psi(\mathbb R^n)$ is closed iff $r(M) = $rank$\,M.$
\item If $P \in \mathcal L_{m,n},$ $M \in \mathbb R^{m \times n},$ $r(M) = n,$ then $F := P \circ \rho_m \circ M \in \mathcal T_n$ and there exist $B \in GL(n, \mathbb R)$ and $Q \in \mathcal L_{n,n}$ such that $F \circ B  = Q \circ \rho_n.$ Hence if $\mathcal N(F)$ is unfolded then $\mathcal N(Q)$ is unfolded, $0 \leq |\Lambda(Q)| < \infty,$ and 
\begin{equation}\label{eqn1.4}
	\Lambda(F) = \bigcup_{\ell = 1}^L \, (B\mathbb Z^n + B\mu_\ell)
\end{equation}
where $L := |\Lambda(Q)|,$ $\{\mu_1,...,\mu_L\} \subset \mathbb C^n$ with $\{\rho_n(\mu_1),...,\rho_n(\mu_L\} = \Lambda(Q),$ and if $\lambda \in B\mathbb Z^n + B\mu_\ell$ then its multiplicity $m_\lambda = m_{\rho_n(\mu_\ell)}.$
\item If $F \in \mathcal T_n$ and $m := $rank$\,G(F)$ 
then there exists $P \in \mathcal L_{m,n},$ and $M \in \mathbb R^{m \times n}$ whose rows are linearly independent over $\mathbb Q$ such that 
$F := P \circ \psi$ where
$\psi := \rho_m \circ M.$ 
\end{enumerate} 
\end{prop}
Proof. 1. Proved in (\cite{lawton}, Lemma 1). 
\\
2. Follows since Ker$(\psi)$ is an abelian Lie group (since it is a closed subgroup of the Lie group $\mathbb R^n$) and $\psi(\mathbb R^n) \simeq \mathbb R^n/\hbox{Ker}\,\psi.$
\\ 
3. Choose $B \in B \in GL(n,\mathbb R)$ such that the columns of
$N := MB \in \mathbb Z^{m \times n}$ form a basis for the rank $n$ $\mathbb Z$--module $Z^m \cap M\mathbb R^n.$
Then $F \circ B = P \circ \rho_m \circ N = P \circ \widetilde N \circ \rho_n = Q \circ \rho_n$ where $Q := P \circ \widetilde N.$ If $\mathcal N(F)$ is unfolded Proposition \ref{prop1.1} implies that $N(Q)$ is unfolded $,\Lambda(Q)$ is finite, and a computation gives (\ref{eqn1.4}).
\\
4. $G(F)$ is a finitely-generated subgroup of $\mathbb R^n$ so it has a $\mathbb Z$--basis 
$\{g_1,...,g_m\}$ by the fundamental theorem of abelian groups. This classic result, proved directly in (\cite{stillwell}, 5.2.5), also follows from the group presentation derived in 1861 by Smith \cite{smith}, (\cite{newman}, Theorem II.p).  Define $M^T := [g_1,...,g_m]^T.$ Then $G(F) = M^T\mathbb Z^m$ 
hence the columns of $M^T$ and rows of $M$ are independent over $\mathbb Q.$ Each component $f_j$ of $F$ is a linear combination of monomials
$h_\omega(z) := e^{2\pi i <\omega,z>}, \omega \in \Omega(F).$ Every $\omega \in \Omega(F)$ has a unique representation 
$\omega = M^Tk, k \in \mathbb Z^m$ since the columns of 
$M^T$ are linearly independent over $\mathbb Q.$ Since $<\omega,z> = <k,Mz>,$ $h_\omega = p_k \circ \rho_m \circ M$ where $p_k \in \mathcal L_{m,1}$ is $p_k(z) := z^k, z \in \mathbb C^{*m}.$ 
\begin{remark}\label{rem1.1}
$(\mathbb T^m,\psi)$ is a compactification of $\mathbb R^n$ 
and $P \circ \psi$ is a representation of the Bohr (or uniformly) almost periodic function $F$ \cite{bohr}.
\end{remark}
Below $L, K, K_1,...,K_n$ are compact convex subsets of $\mathbb R^n,$ $+$ is Minkowski sum, and $V_n$ is $n$-dimensional volume, and 
$\lambda, \mu, \lambda_1,...,\lambda_n \geq 0.$ 
\begin{prop}\label{prop1.3}
$V_n(\lambda_1K_1+\cdots \lambda_nK_n)$ is a homogeneous polynomial
of degree $n$ in $\lambda_1,...,\lambda_n.$ 
\end{prop}
Proof. Minkowski \cite{minkowski1,minkowski2}. 
Ewald (\cite{ewald}, p. 116).
\begin{defi}\label{def1.2} 
The mixed volume $V(K_1,...,K_n) := \frac{1}{n!} \times$ coefficient of $\lambda_1\cdots \lambda_n$
in $V_n(\lambda_1K_1+\cdots +\lambda_nK_n).$
\end{defi}
\begin{prop}\label{prop1.4}
The mixed volume satisfies:
\begin{enumerate}
\item $V(K,...,K) = V_n(K).$
\item If $A = (a_{i,j}) \in \mathbb R^{n \times n}$ and $K_1,...,K_n$ are rectangular bodies parallel to the coordinate axess and
the length of $K_i$ along the $j$-th axes equals $a_{i,j}$ then $n! V(K_1,...,K_n) = $ permanent$\,A.$ 
\item Polarization Identity: $n! \, V(K_1,...,K_n) := (-1)^{n-1} \sum_i V_n(K_j) \\
+ (-1)^{n-2} \sum_{i < j} V_n(K_i + K_j) + \cdots + V_n(K_1 + \cdots + K_n)$
\item $V(K_1,...,K_n) \geq 0.$ 
$V(K_1,...,K_n) > 0$ iff there exists points $a_i, b_i \in K_i, i = 1,...,n$ with $\{b_i-a_i: i = 1,...,n\}$ is linearly independent. Then
$V(K_1,...,K_n) \geq V([a_1,b_1],...,[a_n,b_n]) = \frac{1}{n!}\, |\det [b_1-a_1, ...,b_n-a_n] |.$  
\item Minkowski Linear: $V(\lambda K + \mu L,K_2,...,K_n) = \\
\lambda V(K,K_2,...,K_n)+\mu V(L,K_2,...,K_n).$
\item If $K_1,...,K_n$ are polytopes with vertices in $\mathbb Z^n,$ then
$n! \, V(K_1,...,K_n) \in \mathbb Z.$
\end{enumerate}
\end{prop}
Proof. 1. follows since $V_n(\lambda_1K + \cdots + \lambda_nK) = (\lambda_1+\cdots + \lambda_n)^n V_n(K).$ 
\\
2. follows since $S := \lambda_1K_1+\cdots+\lambda_nK_n$ is a rectangular body whose lenght along the $j$ axis equals
$\ell_j := \lambda_1a_{1,j} + \cdots+\lambda_na_{n,j}$
and $V_n(S) = \ell_1\cdots\ell_n.$
\\
Schneider proved 3, 4, 5  
(\cite{schneider}, Lemma 5.1.4, Theorem 5.1.8, Equation 5.26).
Ewald proved 6 (\cite{ewald}, Theorem 3.9, p. 120). 
\begin{prop}\label{prop1.5}
If $Q \in \mathcal L_{n,n}$ and $\mathcal N(Q)$ is unfolded then 
\begin{equation}\label{eqn1.5}
\sum_{\zeta \in \Lambda(Q)} m_\zeta = n!\, V(\mathcal N(Q)).
\end{equation}
\end{prop}
Proved independently by Bernshtein \cite{bernshtein} and Kouchnirenko \cite{kouchnirenko}.
%
%
\section{Density of Zeros}\label{sec2}
Zeros of univariate entire functions were studies for centuries \cite{vanvleck}. 
This section describes zeros of systems of multivariate trigonometric polynomials.
Note that Proposition \ref{prop2.1} is a special case of Proposition \ref{prop2.2} and its proof is much easier.
\begin{prop}\label{prop2.1}
If $F \in \mathcal T_n$ is periodic, $\mathcal N(F)$ is unfolded and
$V(\mathcal N(F)) > 0,$ then $\Im(\Lambda)$ is bounded and
$\Lambda(F)$ is nonempty, discrete, and its density
\begin{equation}\label{eqn2.1} 
 d(\Lambda(F)) := \lim_{R \rightarrow \infty} 
 \frac{1}{V_n(B(0,R) \cap \mathbb R^n)} 
 \sum_{\lambda \in \Lambda(F),\, ||\lambda||_2 < R} m_\lambda = n!\,V(\mathcal N(F)).
\end{equation}
\end{prop}
Proof. Propositions \ref{prop1.2} and
\ref{prop1.5} give 
$d(\Lambda(F)) = |\det B|^{-1}n! V(\mathcal N(Q)).$ The conclusion follows since $V(\mathcal N(Q))=V(B^T\mathcal N(f_1),...,B^T\mathcal N(f_n)) = |B|\, V(\mathcal N(F).$
\begin{prop}\label{prop2.2}
If $F \in \mathcal T_n,$ $\mathcal N(F)$ is unfolded and
$V(\mathcal N(F)) > 0,$ then $\Im(\Lambda)$ is bounded and
$\Lambda(F)$ is nonempty, discrete, and $d(\Lambda(F)) = n!\, V(\mathcal N(F)).$
\end{prop}
Proof. Follows from formuli of Gelfond  \cite{gelfond1} and Kazarnovskii \cite{kazarnovskii1, kazarnovskii2} for the density  of zeros of systems of holomorphic almost periodic functions.
\\ \\
$\mathcal A(P) := \{ (\ln |z_1|,...\ln |z_m|) : z \in \Lambda(P)\}$
is the amoeba of $P.$ Amoebas were introduced by Gelfand, Kapranov and Zelevinsky \cite{gelfand}.
\begin{defi}\label{defi3.1}
If $M \in \mathbb R^{m \times n}$ then $P \in \mathcal L_{m,n}$ is $M$--stable if 
\begin{equation}\label{eqn4.3}
	\mathcal A(P) \cap M\mathbb R^n = \{0\}.
\end{equation}
\end{defi}
Clearly $P$ is $M$--stable iff $\Lambda(P \circ \rho_m \circ M) \subset \mathbb R^n.$
The stable pairs of polynomials used by Kurasov and Sarnak \cite{kurasovsarnak} correspond to $M$--stable where $n = 1$ and the entries of $M$ are nonzero with the same sign.
\begin{prop}\label{prop2.3}
Nonreal roots of $F \in \mathcal T_n$ are empty or have positive density. All roots of $F$ are in $\mathbb R^n$ iff the density of real roots equal the density of all roots.
\end{prop}
Proof. If $\lambda \in \Lambda(F) \backslash \mathbb R^n$ let $r > 0$ so $B(\lambda,r) \cap \mathbb R^n = \phi.$ The Martinelli-Bochner integral representation (\cite{shabat}, Theorem 1, p. 157) for $m_\lambda$ and almost periodicity of $F$ imply there exist $S \subset \mathbb R^n,$ compact $K \subset \mathbb R^n$ with $S + K = \mathbb R^n$ and $F(s+\lambda) = 0, s \in S.$ Then $S$ and hence $S + \lambda$ have positive density thus proving the first assertion. The second assertion since the density of all zeros is the sum of the density of real zeros, which exists by (\cite{lawton}, Theorem 6), and the density of nonreal zeros.
%
%
\section{Generalized PSF}\label{sec3}
$\mathcal S(\mathbb R^n)$ is the Schwartz space of smooth, i. e. infinitely differentiable, functions all of whose derivatives decay fast \cite{schwartz}. Its dual space $\mathcal S^*(\mathbb R^n)$ is the space of tempered distributions. The Fourier transform is a continuous bijection of $\mathcal S(\mathbb R^n)$ onto itself and extends by duality to a continuous bijection of $\mathcal S^*(\mathbb R^n).$ $\mathcal S_c(\mathbb R^n) \subset \mathcal S(\mathbb R^n)$ is its dense subspace of compactly supported functions. The 
Fourier-Laplace transform $\widehat h : \mathbb C^n \mapsto \mathbb C$ of $h \in \mathcal S_c(\mathbb R^n)$
is
\begin{equation}\label{eqn3.1}
	\widehat h(z) := \int_{x \in \mathbb R^n} h(x)\, e^{-2\pi i \, <x, z>}\, dx, \ \ z \in \mathbb C^n.
\end{equation}
The Paley--Wiener-Schwartz theorem (\cite{hormander}, Theorem 7.3.1) implies that for every
$N > 0$ there exist 
$\gamma > 0$ 
and  
$C_N > 0$ such that
\begin{equation}\label{eqn3.2}
	|\widehat h(z)| \leq C_N(1+|z|)^{-N} e^{\gamma||\Im z||}, \ \ z \in \mathbb C^n,
\end{equation}
and conversely, this condition implies that $h \in \mathcal S_c(\mathcal R^n).$
\begin{prop}\label{prop3.1}
Every $F \in \mathcal T_n$ for which
$\mathcal N(F)$ is unfolded and $V(\mathcal N(F)) > 0$ defines
$\widetilde F \in \mathcal S_c^*(\mathbb R^n)$ by 
\begin{equation}\label{eqn3.3}
	\widetilde F(h) := \sum_{\lambda \in \Lambda(F)} m_\lambda \widehat h(\lambda), \ \ h \in \mathcal S_c(\mathbb R^n).
\end{equation}
\end{prop}
Proof. Proposition \ref{prop1.1} implies $\Im(\Lambda(F))$ is bounded and Proposition \ref{prop2.2} implies $\Lambda(F)$ has finite density, hence (\ref{eqn3.2}) implies that (\ref{eqn3.3}) converges absolutely.
\\ \\
Clearly $F \in \mathcal T_n$ is periodic iff its period group $\Gamma(F) = B\mathbb Z^n$
for $B \in GL(n,\mathbb R)$ iff $F = Q \circ \rho_n \circ B^{-1}$ where $Q = (q_1,...,q_n) \in \mathcal L_{n,n}.$ Propositions \ref{prop3.2} and \ref{prop3.3} compute $\widetilde F(h)$ under the hypothesis that $F$ has this representation where $\mathcal N(Q)$ and hence $\mathcal N(F)$ are unfolded and their mixed volumes are positive. 
\begin{prop}\label{prop3.2}
Under the preceding hypotheses
\begin{equation}\label{eqn3.4} 
\widetilde F(h) =  |\det B|^{-1}\sum_{\zeta \in \Lambda(Q)} m_\zeta\, R_h(\zeta)
\end{equation}
where $m_\zeta$ is the multiplicity of $\zeta,$
and $R_h \in \mathcal L_{n,1}$ is
\begin{equation}\label{eqn3.5}
R_h(z) :=  \sum_{k \in \mathbb Z^n} 
	h(B^{-T}k)\, z^{-k}, \ \ z \in \mathbb C^{*n}.
\end{equation}
\end{prop}
Proof Proposition \ref{prop1.2} implies that 
$\widetilde F(h) = \sum_\lambda m_\lambda \widehat h(\lambda)$ where 
$$\lambda \in \bigcup_{\ell = 1}^L \, (B\mathbb Z^n + B\mu_\ell)$$ 
where $L = |\Lambda(Q)|$ and 
$\Lambda(Q) = \{\rho_n(\mu_\ell), \ell = 1,...,L\}.$
Since $F$ is invariant under translation by elements in $B\mathbb Z^n$ 
and $\rho_n\circ B^{-1}$ is locally a holomorphic homeomorphism,
$m_\lambda = m_{B\mu_\ell} = m_{\rho_n(\mu_\ell)}$ 
for $\lambda \in \mathbb BZ^n + B\mu_\ell,$ $\ell = 1,...,L.$ Therefore
\begin{equation}\label{eqn3.6}
\widetilde F(h) = \sum_{\ell = 1}^L m_{\rho_n(\mu_\ell)} 
\sum_{k \in \mathbb Z^n} \widehat h(Bk+B\mu_\ell).
\end{equation}
The proof is finished since the classical PSF (\cite{hormander}, Theorem 7.2.1) gives
\begin{equation}\label{eqn3.7} 
	\sum_{k \in \mathbb Z^n} \widehat h(Bk+B\mu_\ell) =
|\det B|^{-1} \sum_{k \in \mathbb Z^n} 
	h(B^{-T}k)e^{-2\pi i \, k \cdot \mu_\ell}, \ \ \ell = 1,...,L.
\end{equation}
Proposition \ref{prop3.2} implies  that
\begin{equation}\label{eqn3.8} 
		\widetilde F(h) = \sum_{\zeta \in \Lambda(Q)} (2\pi i)^{-n}\int_{G_\zeta} \omega
\end{equation}
where the Grothendieck $n$--cycle at $\zeta,$
\begin{equation}\label{eqn3.9}
	G_\zeta := \{z \in \mathbb C^{*n} : ||z-\zeta||_2 \leq \epsilon, |q_j(z)| = \epsilon_j, j = 1,...,n\},
\end{equation}
is oriented so $d(\arg q_1) \wedge \cdots \wedge d(\arg q_n) > 0$ on
$G_\zeta,$
and the $n$--form
\begin{equation}\label{eqn3.10}
	\omega := R_h |\det B|^{-1} \frac{dq_1 \wedge \cdots \wedge dq_n}{q_1\cdots q_n}.
\end{equation}
Direct computation gives $\omega = H \frac{dz_1}{z_1} \wedge \cdots \wedge \frac{dz_n}{z_n}$ where
\begin{equation}\label{eqn3.11}
	H :=  
	R_h \frac{
	|\det B|^{-1} \, z_1 \cdots z_n \det \left(\partial_{z_j} q_k\right)}
	{q_1\cdots q_n}.
\end{equation}
Observe that 
\begin{equation}\label{eqn3.12}
\mathcal N(q_1\cdots q_n) = 
\mathcal N(q_1) + \cdots + \mathcal N(q_n).
\end{equation}
For every $v \in \mathcal V(q_1\cdots q_n),$ Gelfond and Khovanskii (\cite{gelfondkhovanskii1}, 1.9) construct an $n$--cycle disjoint from $N_v \subset \mathbb C^{*n}\backslash \Lambda(q_1\cdots q_n)$ of the form $N_v :=  c(v)\mathbb T^n,$ where $c(v) 
\in \mathbb C^{*n},$ and oriented so $\frac{dz_1}{z_1} \wedge \cdots \wedge \frac{dz_n}{z_n} > 0$ on $N_v.$ They define the residue 
\begin{equation}\label{eqn3.13}
	\hbox{res}_v \, \omega := (2\pi i)^{-n}\int_{N_v} \omega.
\end{equation}
Clearly 
\begin{equation}\label{eqn3.14}
	\hbox{res}_v \omega = H_v
\end{equation}
where we define $H_v$ to be the constant term of the Laurent expansion of $H$ at $v.$
Gelfond and Khovanskii (\cite{gelfondkhovanskii1}, 1.10) prove that
\begin{equation}\label{eqn3.15}
	\sum_{\zeta \in \Lambda(Q)} G_\zeta
	\sim
	(-1)^n \sum_{v \in \mathcal V(q_1\cdots q_n)} k_v\, N_v.
\end{equation}
where $\sim$ means homologous in $\mathbb C^{*n} \backslash \Lambda(q_1\cdots q_n),$ and $k_v \in \mathbb Z$ is the 
combinatorial coefficient \cite{gelfond2} of $\mathcal N(Q)$ at $v.$ Therefore Stoke's theorem, (\ref{eqn3.8}), (\ref{eqn3.13}), (\ref{eqn3.14}) and (\ref{eqn3.15}) give
\begin{equation}\label{eqn3.16}
		\widetilde F(h) = (-1)^n \sum_{v \in \mathcal V(q_1\cdots q_n)} k_v \, H_v.
\end{equation}
We observe that $H_v = (H \circ \rho_n \circ B^{-1})_{w},$
the constant term in the Fourier expansion of $H \circ \rho_n \circ B^{-1}$ at $w := B^{-T}v \in \mathcal V(f_1\cdots f_n).$ Moreover the chain rule for differentiation gives
\begin{equation}\label{eqn3.17}
		H \circ \rho_n \circ B^{-1} = 
		(R_h\circ \rho_n \circ B^{-1})\, 
		\frac{(2\pi i)^{-n} \det \left(\partial_{x_i} f_j\right)}{f_1\cdots f_n}. 
\end{equation}
Exressing $f_1 \cdots f_n(x) = ce^{2\pi i w \cdot x}(1 - g(x))$
gives 
\begin{equation}\label{eqn3.18}
		\frac{1}{(f_1\cdots f_n)(x)} = c^{-1}e^{-2\pi i w\cdot x}
		\left( 1 + \sum_{k=1}^\infty g(x)^k \right), \ \ x \in B\rho_n^{-1}(N_v).
\end{equation}
Since 
$\Omega(g) = (\Omega(f_1\cdots f_n) - \{w\})\backslash \{0\}$ 
is a finite subset of a pointed cone with $0$ removed, there exists a discrete $S_w \subset \mathbb R^n,$ $a_w : S_w \mapsto \mathbb C,$ $N > 0$ and $C > 0$ satisfying 
\begin{equation}\label{eqn3.19}
\sum_{s \in S_w \cap B(0,R)} |a_w(s)| < C\, (1+R)^N, \ \ R > 0
\end{equation}
such that
\begin{equation}\label{eqn3.19}
\frac{(2\pi i)^{-n} \det \left(\partial_{x_i} f_j\right)}{f_1\cdots f_n}(x)
= \sum_{s \in S_w} a_w(s) \, e^{2\pi i s \cdot x}, \ \  x \in B\rho_n^{-1}(N_v).
\end{equation}
Therefore (\ref{eqn3.14}) and (\ref{eqn3.16})-(\ref{eqn3.20}) give
\begin{equation}\label{eqn3.20}
H_v = \sum_{s \in S_w} a_w(s) h(s).
\end{equation}
\begin{prop}\label{prop3.3} Under the hypotheses in Proposition \ref{prop3.2}
there exists a discrete $S \subset \mathbb R^n$ and $a : S \mapsto \mathbb C$ such that
\begin{enumerate}
\item $\sum_{s \in S} |a(s)|\, \delta_s \in \mathcal S^*(\mathbb R^n),$ 
\item $\widetilde F(h) = \sum_{s \in S} a(s)\, h(s).$
\end{enumerate} 
\end{prop}
Proof. Let $S = \bigcup_{w \in \mathcal V(f_1\cdots f_n)} S_w$
and $a = \sum_{w \in \mathcal V(f_1\cdots f_n)} k_w\, a_w.$
Here we observe that $k_w = k_v$ where $w = B^{-T}v$ and
$k_w$ depends on the combinatorics of 
$\mathcal N(f_1\cdots f_n) = 
\mathcal N(f_1) + \cdots +\mathcal N(f_n).$
Then (\ref{eqn3.19}) implies 1 and 
(\ref{eqn3.21}) implie 2.
%
%
\begin{theo}\label{thm3.1}
If $F \in \mathcal T_n,$ $\mathcal N(F)$ is unfolded, $V(\mathcal N(F)) > 0,$ and $h \in \mathcal S_c$ then there exist a discrete $S \subset \mathbb R^n$ and
$a : S \rightarrow \mathcal C$ such that
\begin{enumerate}
\item $\zeta := \sum_{s \in S} a(s) \delta_s \in \mathcal S^*(\mathbb R^n)$ and
$\widetilde F(h) = \zeta(h),$
\item $|\zeta| := \sum_{s \in S} |a(s)| \delta_s
\in  \mathcal S^*(\mathbb R^n).$
\end{enumerate}
\end{theo}
Proof. Represent $F = P \circ \rho_m \circ M$ as in part 4 of Proposition \ref{prop1.2},
let $M_k \in \mathbb Q^{m \times n}$ be a sequence such that $||M-M_k||_2 \rightarrow 0,$
and define $F_k := P \circ \rho_m \circ M.$ Then (i) for sufficiently large $k,$ $F_k$ satisfies the hypotheses in Proposition \ref{prop3.3} and (ii) $F_k \rightarrow F$ uniformly on compact subsets of $\mathbb C^n.$ (i), (ii), and Proposition \ref{prop3.3} imply that for sufficiently large $k$ there exist discrete $S_k \subset \mathbb R^n$ and $a_k : S_k \mapsto \mathbb C$ such that
$\zeta_k := \sum_{s \in S_k} a_k(s) \delta_s \in \mathcal S^*(\mathbb R^n)$ and 
\begin{equation}\label{eqn3.21}
\widetilde {F_k}(h) = \zeta_k(h), \ \ h \in \mathcal S_c(\mathbb R^n).
\end{equation}
and $|\zeta_k| \in S^*(\mathbb R^n).$ Property 2 implies
$S_k$ converges pointwise to a discrete $S \subset \mathbb R^n$ and $a_k$ converges pointwise to $a : S \mapsto \mathbb C.$ Define 
$\zeta := \sum_{s \in S} a(s) \delta_s \in S^*(\mathbb R^n).$ Then
\begin{equation}\label{eqn3.22}
\widetilde F(h) = \lim_{k \rightarrow \infty} \widetilde {F_k}(h) = \zeta(h), \ \ h \in \mathcal S_c(\mathbb R^n),
\end{equation}
and $|\zeta| \in S^*(\mathbb R^n).$ This concludes the proof.
%
%
\section{Fourier Quasicrystals}\label{sec4}
%
%
\begin{theo}\label{thm4.1}
If $F \in \mathcal T_n,$ $\mathcal N(F)$ is unfolded, $V(\mathcal N(F)) > 0$
and $\Lambda(F) \subset \mathbb R^n,$ then
$\mu := \sum_{\lambda \in \Lambda(F)} m_\lambda \, \delta_\lambda,$ 
where $m_\lambda \in \mathbb N$ is the multiplicity of $\lambda,$ is a FQ. 
\end{theo}
Proof. Let $\zeta := \sum_{s \in S} a(s) \delta_s \in \mathcal S^*(\mathbb R^n)$
be the atomic measure in Theorem \ref{thm3.1} such that $\zeta, |\zeta| \in \mathcal S^*(\mathbb R^n)$ and for every $h \in \mathcal S_c(\mathbb R^n)$
\begin{equation}\label{eqn4.1}
\widehat \mu(h) = \mu(\widehat h) = \widetilde F(h) = \zeta(h).
\end{equation}
Since $\Lambda(F) \subset \mathbb R^n$ and $\mathcal S_c(\mathbb R^n)$ is dense in $\mathcal S(\mathbb R^n),$ (\ref{eqn4.1}) holds for all $h \in  \mathcal S(\mathbb R^n)$ hence $\mu$ is a FQ.  
%
%
\begin{example}\label{ex4.1}
Let $n \geq 2, m = n+1,$ $s_1,...,s_n \in (-1,1)\backslash \{0\},$
and $P := (p_1,...,p_n) \in \mathcal L_{m,n}$ where
\begin{equation}\label{eqn4.1}
	p_j(z_1,...,z_m) = z_j(1+s_jz_m)-z_m - s_j, \ \ j = 1,...,n,
\end{equation}
$b_j < 0$ satisfy $1, b_1, b_2$ are rationally independent and define
$M := \left[
\begin{array}{c}
I_n \\
b^T 
\end{array}
\right],$
and define $\psi := \rho_n \circ M$ and 
$F := P \circ \rho_m \circ M.$ Then 
\begin{enumerate}
\item $P$ is $M$--stable, 
\item $\mu := \sum_{\lambda \in \Lambda(F)} \delta_\lambda$
is a FQ, 
\item $\Lambda(F)$ contains no infinite arithmetic progression. 
\end{enumerate}
\end{example}
Proof of 1.
Since $p_j = 0$ iff $z_j = \frac{z_m - s_j}{1+s_jz_m},$ if $P = 0$ then the signs of 
$\log|z_j|, j =1,...,m$ must be the same. If $x, y \in \mathbb R^n$ 
and $z := \rho_m \circ M(x+iy)$ then $sgn \log |z_j | = sgn (-2\pi (My)_j = -sgn (My)_j, j = 1,...,m.$ There exists a nonzero vector $y$ such that these quantities
have the same sign or equal zero iff the cone spanned by the rows of $M$ lie in a closed half--space. 
\\
Proof of 2.
Clearly $\mathcal N(F) = (M^T \mathcal NM(p_1),...,M^T \mathcal N(p_n)$ is unfolded
and its mixed volume is positive so Theorem \ref{thm3.1} concludes the proof.
\\
Proof of 3.
Since $1,b_1,...,b_n$ are rationally independent $\psi$ is one-to-one.
If $\Lambda(F)$ contains an arithmetic progression $\Gamma_1$ then $\Gamma_1$ is contained in a finite rank
subroup $\Gamma_2$ of $\mathbb R^n.$ Define 
$G_k := \rho_n \circ M(\Gamma_k), k = 1,2.$
Then $G_1 \subset G_2$ and $G_2$ is a finite rank subgroup of $\mathbb T^m.$
Lang'sconjecture (\cite{evertse}, Theorem 10.10.1) implies $G_2 \cap \Lambda(P)$ is 
contained in the finite union of translates of torus subgroups of $\Lambda(P).$ Since 
$\{1\}$ is the only torus subgroup of $\Lambda(P),$ $G_2 \cap \Lambda(P)$ is finite 
hence $G_1$ is finite since $G_1 \subset G_2 \cap \Lambda(P).$ Since $\rho_m \circ M$ 
is one-to-one, $\Gamma_1$ is finite.
%
%
\begin{example}\label{ex4.2}
Let $n \geq 2, m = n+1,$ $\delta \in (-1,1)\backslash \{0\},$
$P := (p_1,...,p_n) \in \mathcal L_{3,2}$
\begin{equation}\label{eqn4.2}
	p_j(z_1,...,z_m) := z_j  - z_j^{-1} - \delta (z_3 -z_3^{-1}), \ \ j = 1,...,n.
\end{equation}
Let $e_1,...,e_n$ be the standard basis vectors for $\mathbb R^n.$
Let $b \in \mathbb R^n$ so $\sum_{j=1}^n |b_j|< 1$ and $1, b_1,...,b_n$ are rationally independent.Define $M \in \mathbb R^{m \times n}$ so $M^T := [I_n \, b]$ and $\psi := \rho_n \circ M$ and $F := (f_1,...,f_n) := P \circ \psi.$ 
Then
\begin{enumerate}
\item $\psi (R^n)$ is dense in $\mathbb T^n.$ 
\item $\psi (R^n)$ is transversal to $\Lambda(P) \cap \mathbb T^n.$
\item The density of $\Lambda(F) \cap \mathbb R^n$ equals
$2^n.$
\item The density of $\Lambda(F)$ equals $2^n.$
\item $P$ is $M$--stable. 
\item $\mu := \sum_{\lambda \in \Lambda(F)} \delta_\lambda$
is a FQ. 
\item $\Lambda(F)$ contains no infinite arithmetic progression. 
\end{enumerate}
\end{example}
1. follows from Proposition \ref{prop1.1} since each $b_j$ being irrational implies that the rows of $M$ are independent over 
$\mathbb Q.$  \\
2. Parameterize $z_j = e^{i \theta_j}, j = 1,...,m.$ $\Lambda(P) \cap \mathbb T^n$ is the union of $2^n$ 
loops defined by equations $\sin \theta_j = \pm \delta \sin \theta_m, j  = 1,...,n$ each homotopic to the circle subgroup $C = \{z \in \mathbb T^n:z_j = 1, j = 1,...,n\}.$ The vectors
$v = \left[\frac{d\theta_1}{d\theta_m},...,\frac{d\theta_n}{d\theta_m}, 1 \right]^T$ are tangents to these curves where
$\frac{d\theta_j}{d\theta_m} = 
\pm \delta \frac{1-\sin^2 \theta_j}{1-\delta^2 \sin^2 \theta_j}, j = 1,...,n$ and all of the first $n$ entries of this vector have modulus $M 1.$ The vector $[1 -b^T]^T$ is normal to $\psi(\mathbb R^n)$ and cannot be normal to the vectors $v$ hence $\psi (R^n)$ is transversal to $\Lambda(P) \cap \mathbb T^n.$ \\
3. follows from the above two properties by (\cite{lawton}, Theorem 6). \\
4. $\mathcal N(f_j) = \hbox{convex hull }\{\pm e_j, \pm b\}$ which up to translation equals the Minkowski sum $[0, e_j+b] + [0,e_j-b].$
Therefore $\mathcal N(F)$ is unfolded and Proposition \ref{prop1.3} implies that its mixed volume 
\begin{equation}\label{eqn4.3}
\begin{array}{ccc}
	V(\mathcal N(F)) & = & \sum_{s_1,...,s_n \in \{1,-1\}} V([0, e_1+s_1b],...,[0, e_n+s_nb]) \\
	& = & \sum_{s_1,...,s_n \in \{1,-1\}} |\det [e_1+s_1b,..., e_n+s_nb]\, | \\
	& = & \sum_{s_1,...,s_n \in \{1,-1\}} \det [e_1+s_1b,..., e_n+s_nb]\, \\
	& = & \sum_{s_1,...,s_n \in \{1,-1\}} \det [e_1,..., e_n] \\
	& = & \sum_{s_1,...,s_n \in \{1,-1\}} 1 \\
	& = & 2^n. \\
\end{array}
\end{equation}
Since the matrices are real, strictly diagonally dominant with positive diagonal entries, Gershgorin's circle theorem 
(\cite{hornjohnson}, Theorem 6.1.10) implies their eigenvalues have positive real part so their determnants are positive which gives the third equality. The forth equality follows since the determinant is a linear function of each column vector. \\
5. then follows from Proposition \ref{prop2.3} \\
6. follows from Theorem \ref{thm4.1} \\
7. Since $1,b_1,...,b_n$ are rationally independent $\psi$ is one-to-one, the conclusion follows as in Example \ref{ex4.1}
since $\Lambda(P) \cap \mathbb T^n$ does not contain a translated of a circle subgroup.
%
%
\section{Research Questions}\label{sec5}
We raise questions for future research and suggest approaches to study them.
\begin{question}\label{ques5.1}
If $P \in \mathcal L_{m,n},$ $M \in \mathbb R^{m \times n},$ and $P$ is $M$--stable, under what conditions is it $M_1$--stable for $M_1$ sufficiently close to $M?$ 
\end{question}
For $n = 1$ and $m = 2$ every connected component
of $\mathbb R^2 \backslash \mathcal A(P)$ corresponds to a point in
$\mathcal N(P)$ \cite{forsbergpassaretsikh, viro}. If $P$ is $M$--stable then it is $M_1$-stable if $M_1$ succifiently close to $M$ iff each rays in $M\mathbb R \backslash \{0\}$ belongs to the connected component of $\mathbb R^2 \backslash \mathcal A(P)$ corresponding to a vertex of $\mathcal N(P).$ This follows since the other components are either bounded or haves boundaries that are asymptotically parallel. For $n \geq 2$ and $m > n$ Question \ref{ques5.1} is more difficult. However, we suggest that the results of Henriques \cite{henriques} and Bushueva and Tsikh \cite{bushuevatsikh}, which derive analogue of convexity for complements of amoebas of higher dimension, may provide answers.
\begin{question}\label{ques5.2}
Are the sufficient conditions in Theorem \ref{thm4.1} also necessary? 
\end{question}
For $n = 1$ the answer was proved yes by Olevskii and Ulanovskii \cite{olevskiiulanovskii2}. But their proof uses several methods, such as the Weierstrass factorization for univariate entire functions, that do not 
have analogues for multivariate functions. Favorov \cite{favorov1, favorov2} derived conditions for a Bohr almost periodic set 
$\Lambda \subset \mathbb R^n$ to be the zero set of a 
holomorphic system. This fact suggests that
the compactification of $\Lambda,$ i.e. the closure $K$ of 
$\psi(\Lambda)$ in an appropriate compact group $G$ might 
be contained in the zero set of a system $H$ of analytic functions 
on $G.$ The spectrum $S \subset \mathbb R^n$ of $\mu$ is 
discrete but $S$ equals a certain projection of the support of 
the Fourier transform $\widehat \nu$ of a measure $\nu$ 
supported on $K,$ as was shown when $G$ is a torus group 
in \cite{lawton}. This implies that the support of 
$\widehat \nu$ is an extremely sparse subset of the Pontryagin dual 
$\widehat G.$ Since the product $H \nu = 0$ the convolution 
$\widehat \nu * \widehat H = 0.$ Since $H$ is analytic 
$\widehat H$ decays exponentially fast. This suggests 
that $\widehat H$ has finite support so $H$ is a system 
of trigonometric polynomial on $G$ so $G$ can be replaced 
by a torus and $K$ is an algebraic variety and and $\Lambda$ 
is a Bohr almost periodic set of toral type.   
\\ \\
{\bf Acknowledgments} The authors thanks Sergey Favorov, Alexander Olevskii, Yves Meyer, and Peter Sarnak for sharing their knowledge of crystaline measures, Fourier quasicrystals and Lang's conjecture. 

\end{LARGE}
\Finishall
\end{document}